\documentclass{amsart}
\usepackage{amssymb,amsmath,euscript,mathrsfs,bbm,txfonts,verbatim}
\newtheorem{theorem}{Theorem}
\newtheorem{lemma}{Lemma}
\newtheorem{corollary}{Corollary}

\title{The essential norm of a composition operator on the minimal M\"obius invariant space}
\author{Themis Mitsis}
\author{Michael Papadimitrakis}
\address{Department of Mathematics, University of Crete, Knossos Ave., 71409 Iraklio, Greece}
\email{themis.mitsis@gmail.com}
\email{papadim@math.uoc.gr}
\begin{document}
\begin{abstract}
We derive a formula for the essential norm of a composition operator on the
minimal M\"obius invariant space of analytic functions. As an application, we show that
the essential norm of a non-compact composition operator is at least 1. We also obtain
lower bounds depending on the behavior of the symbol near the boundary, and calculate
the order of magnitude of the essential norm of composition operators induced by finite
Blaschke products.
 \end{abstract}

\pagestyle{plain}
\maketitle
Let $\mathbbm D\subset\mathbbm C$ be the open unit disk. For $\alpha\in\overline{\mathbbm D}$, we put
\[\varphi_\alpha(z)=\frac{\alpha-z}{1-\overline \alpha z},\quad z\in\mathbbm D.\]
The minimal space $\mathscr M$ (or analytic Besov-$1$ space) is defined to be the set of all
analytic functions $f$ on $\mathbbm D$ for which there exist a sequence of points $\alpha_n\in\overline{\mathbbm D}$, and a sequence of complex numbers $\lambda=(\lambda(n))_{n=1}^\infty\in \ell^1$ such that
\begin{equation*}
%\label{atomic}
f(z)=\sum_{n=1}^\infty \lambda(n)\,\varphi_{a_n}(z).
\end{equation*}
So, a function in $\mathscr M$ has an ``atomic'' decomposition as a sum of M\"obius transformations. We norm $\mathscr M$ by
\[\|f\|_\mathscr M=\inf\left\{\|\lambda\|_{\ell^1}:f=\sum_{n=1}^\infty \lambda(n)\,\varphi_{\alpha_n},\ \text{for some $\lambda\in \ell^1$, $\alpha_n\in\overline{\mathbbm D}$}\right\}.\]

The minimal space was introduced and extensively studied in \cite{afp}, where it was shown that if one defines appropriately the notion of a ``M\"obius invariant space'', then $\mathscr M$ is the smallest one. In fact, its norm is stronger than the norm of any other such space. Note that functions in $\mathscr M$ extend continuously to the boundary, so $\mathscr M$ is a ``boundary regular'' space. Moreover, as shown in \cite{afp}, $\mathscr M$ coincides with the set of all analytic functions on $\mathbbm D$ with integrable second derivative. More specifically, there exists a constant $C>0$ such that for every $f\in\mathscr M$
\begin{equation*}
\label{pe1}
C^{-1}\,\|f''\|_1\leq\|f-f(0)-f'(0)\,z\|_\mathscr M\leq C\,\|f''\|_1,
\end{equation*}
where
\[\|f''\|_1=\iint\limits_\mathbbm D|f''(z)|\, dA(z),\]
and $dA$ is normalized area measure. Every $f\in\mathscr M$ can be recovered from its second derivative by means of the  formula
\[f(z)=f(0)+f'(0)\, z-\iint\limits_\mathbbm D\frac{\varphi_\alpha(z)}{\overline\alpha}f''(\alpha)\, dA(\alpha).\] 
In what follows, for positive $x,y$, the symbols $x\simeq y$ and  $x\lesssim y$ mean  $C^{-1}x\leq y\leq C\, x$ and $x\leq C\, y$ respectively, where $C>0$ is an absolute numerical constant, not necessarily the same each time it occurs.

Now let $\psi$ be an analytic map of $\mathbbm D$ into itself. It is clear that the composition operator
\[C_\psi:\mathscr M\to\mathscr M, \quad C_\psi f=f\circ\psi,\]
is bounded if and only if
\begin{equation}
\label{bounded}
\sup_{|\alpha|<1}\|(\varphi_{\alpha}\circ\psi)''\|_1<\infty.
\end{equation}
 Of course, this is equivalent to \[\sup_{|\alpha|<1}\iint\limits_{\mathbbm D}\left|\frac{1-|\alpha|^2}{(1-\overline \alpha\psi)^2}\psi''+2\overline \alpha\frac{1-|\alpha|^2}{(1-\overline\alpha\psi)^3}(\psi')^2\right|\, dA<+\infty.\] In \cite{afp} the authors prove that this is equivalent to
\[\sup_{|\alpha|<1}\iint\limits_{\mathbbm D}\frac{1-|\alpha|^2}{|1-\overline \alpha\psi|^2}|\psi''|\, dA<+\infty\]
together with
\[\sup_{|\alpha|<1}\iint\limits_{\mathbbm D}\frac{1-|\alpha|^2}{|1-\overline \alpha\psi|^3}|\psi'|^2\, dA<+\infty.\]
One direction of this equivalence is, of course, trivial and the other direction is done in \cite{afp} by proving, via interpolation with the Dirichlet space, that, if $C_\psi$ is bounded on $\mathscr M$, then the second of the last two relations holds and hence also the first. Here is an alternative proof. If $C_\psi$ is bounded on $\mathscr M$, then, since \[|f(0)|+\iint\limits_{\mathbbm D}|f'|^2\, dA\lesssim\|f\|_\mathscr M,\] we get that
\[\sup_{|\alpha|<1}\iint\limits_{\mathbbm D}\frac{(1-|\alpha|^2)^2}{|1-\overline\alpha\psi|^4}|\psi'|^2\, dA<+\infty.\]
This is equivalent to 
\[\sup_{|\alpha|<1}\iint\limits_{\mathbbm D}\frac{(1-|\alpha|^2)^2}{|1-\overline\alpha z|^4}n_\psi(z)\, dA(z)<+\infty,\]
where $n_\psi(z)$ is the cardinality of the inverse image of $z$ under $\psi$. Now, it is standard to show that this is equivalent to 
\[\sup_I\frac 1{|I|^2}\iint\limits_{S(I)}n_\psi(z)\, dA(z)<+\infty,\]
where $I$ ranges over all arcs on the unit circle $\mathbbm T$, $|I|$ is the length of $I$,  and $S(I)$ is the usual Carleson square over $I$. Finally, it is also standard to show that this is equivalent to
\[\sup_{|\alpha|<1}\iint\limits_{\mathbbm D}\frac{1-|\alpha|^2}{|1-\overline\alpha z|^3}n_\psi(z)\, dA(z)<+\infty,\]
which is equivalent to the original \[\sup_{|\alpha|<1}\iint\limits_{\mathbbm D}\frac{1-|\alpha|^2}{|1-\overline\alpha\psi|^3}|\psi'|^2\, dA<+\infty.\] 
A sufficient condition for the boundedness of $C_\psi$, involving the integral means of $\psi''$, has been obtained  by Blasco in \cite{blasco}, extending a series of partial results in \cite{afp}. As far as compactness is concerned, the minimal space satisfies the conditions of theorem 2.1 in \cite{shapiro}, therefore $C_\psi$ is compact if and only if $\|\psi\|_\infty<1$. Moreover, Wulan and Xiong proved in \cite{wulan} that $C_\psi$ is compact if and only if the ``little Oh'' version of (\ref{bounded}) holds, namely
\[\lim_{|\alpha|\to1}\|(\varphi_\alpha\circ\psi)''\|_{1}=0.\]

Here we are interested in estimating the essential norm of $C_\psi$ which is defined to be the distance of $C_\psi$ to the subspace of compact operators, that is
\[\|C_\psi\|_e=\inf\{\|C_\psi-K\|:\text{$K$ compact}\}.\]
We will prove the following asymptotic estimate. The result in \cite{wulan} is, of course, a special case.  

\begin{theorem}
\label{main}
Let $C_\psi:\mathscr M\to\mathscr M$ be a bounded composition operator. Then
\[\|C_\psi\|_e\simeq\limsup_{|\alpha|\to1}\|(\varphi_\alpha\circ\psi)''\|_{1}:=\lim_{s\to1}\,\sup_{|\alpha|>s}\|(\varphi_\alpha\circ\psi)''\|_{1}.\]
\end{theorem}
In the proof we will make use of a simple observation: A  bounded sequence $f_n\in\mathscr M$ converges to zero in the weak-* topology if and only if it converges to zero uniformly on compact sets. To see this, note that by \cite{afp},  the pre-dual of $\mathscr M$ may be identified with the little Bloch space. More specifically, every element $\Lambda$ of the dual of the little Bloch space can be represented as
\[\Lambda(b)=\Lambda_f(b):=\langle b,f\rangle+b(0)\,\overline{f(0)},\]
for some $f\in\mathscr M$. Here
\[\langle b,f\rangle=\iint\limits_\mathbbm Db'\,\overline{f'}\, dA,\]
is the invariant pairing. Now for $a,z\in\mathbbm D$   let
\[b_{a}(z)=\log\frac1{1-\overline az}.\]
Then $b_a$ is a little Bloch function. So, assuming that $f_n$ is weak-* null, we have that $\Lambda_{f_n}(b_a+1)\to0$ for all $a$. However
\[\Lambda_{f_n}(b_a+1)=\langle b_a,f_n\rangle+\overline{f_n(0)},\]
and a calculation shows that $\langle b_a,f_n\rangle=\overline{f_n(a)-f_n(0)}$. Consequently, $f_n\to0$ pointwise. Since $f_n$ is bounded, an application of Cauchy's theorem shows that we actually have uniform convergence on compact sets.  To prove the converse, it is enough to show that $\Lambda_{f_n}(z^k)\to0$ for $k=0,1,2,\dots$, because polynomials are dense in the little Bloch space. However,
\[\Lambda_{f_n}(1)=\overline{f_n(0)},\quad \Lambda_{f_n}(z^k)=\langle z^k,f_n\rangle=\frac{\overline{f_n^{(k)}(0)}}{(k-1)!}.\]

\begin{proof}[Proof of theorem \ref{main}]
For $R\in(0,1)$ we will use the notation
\[L_R=\{z\in\mathbbm D:|\psi(z)|<R\},\quad U_R=\{z\in\mathbbm D:|\psi(z)|\geq R\}.\]

 To prove the lower bound
\begin{equation}
 \label{e1}
\|C_\psi\|_e\gtrsim\limsup_{|\alpha|\to1}\|(\varphi_\alpha\circ\psi)''\|_{1},
\end{equation}
choose a sequence $\alpha_n\in\mathbbm D$ with $|\alpha_n|\to1$, such that
\[\limsup_{|\alpha|\to1}\|(\varphi_\alpha\circ\psi)''\|_{1}=\lim_n\|(\varphi_{\alpha_n}\circ\psi)''\|_{1},\]
and let $K:\mathscr M\to\mathscr M$ be a compact operator. Without loss of generality, we may assume that $\|K\varphi_{\alpha_n}-f\|_\mathscr M\to0$, for some $f\in\mathscr M$. Then 
\begin{equation}
\label{e2}
\|C_\psi-K\|\geq\|C_\psi\varphi_{\alpha_n}-K\varphi_{\alpha_n}\|_\mathscr M\gtrsim\|(\varphi_{\alpha_n}\circ\psi)''-f''\|_{1}-\|K\varphi_{\alpha_n}-f\|_\mathscr M.
\end{equation}
We may further assume that $f''$ is not identically zero, otherwise (\ref{e1}) follows from (\ref{e2}) upon taking the limit as $n\to\infty$. Now let $\varepsilon>0$, and fix $R\in(0,1)$ such that 
\[\iint\limits_{U_R}|f''|\, dA<\varepsilon.\]
On the other hand 
\begin{align*}
 \iint\limits_{L_R}|(\varphi_{\alpha_n}\circ\psi)''|\, dA&=\iint\limits_{L_R}\left|\frac{1-|\alpha_n|^2}{(1-\overline a_n\psi)^2}\psi''+2\overline a_n\frac{1-|\alpha_n|^2}{(1-\overline a_n\psi)^3}(\psi')^2\right|\, dA\\
&\lesssim\frac{1-|\alpha_n|^2}{(1-R)^3}\left(\|\psi''\|_{1}+\|\psi'\|^2_{2}\right).
\end{align*}
Therefore, for  $n$ large enough
\[\iint\limits_{L_R}|(\varphi_{\alpha_n}\circ\psi)''|\, dA<\frac12\iint\limits_{L_R}|f''|\, dA.\]
Consequently
\begin{align*}
 \|(\varphi_{\alpha_n}\circ\psi)''-f''\|_{1}&=\iint\limits_{L_R}|(\varphi_{\alpha_n}\circ\psi)''-f''|\, dA+\iint\limits_{U_R}|(\varphi_{\alpha_n}\circ\psi)''-f''|\, dA\\
&\geq\iint\limits_{L_R}|f''|\, dA-\iint\limits_{L_R}|(\varphi_{\alpha_n}\circ\psi)''|\, dA+\iint\limits_{U_R}|(\varphi_{\alpha_n}\circ\psi)''|\, dA-\varepsilon\\
&>\|(\varphi_{\alpha_n}\circ\psi)''\|_{1}-\varepsilon.
\end{align*}
It follows that
\[\|C_\psi-K\|\gtrsim\limsup_{|\alpha|\to1}\|(\varphi_\alpha\circ\psi)''\|_{1}-\varepsilon.\]
The idea of the preceding argument is that for $\alpha$ close to the boundary, $(\varphi_\alpha\circ\psi)''$ and $f''$, have, loosely speaking of course, ``disjoint supports''.

To prove the upper bound
\begin{equation*}
\|C_\psi\|_e\lesssim\limsup_{|\alpha|\to1}\|(\varphi_\alpha\circ\psi)''\|_{1},
\end{equation*}
let $\varepsilon,s,R,R'\in(0,1)$ with  $R<R'$, and for $f(z)=\sum_{k=0}^\infty a_kz^k\in\mathscr M$ let
\[T_nf(z)=\sum_{k=0}^n\left(1-\frac k{n+1}\right)a_kz^k,\quad Pf(z)=a_0+a_1z.\]
Then $T_nf\to f$ in $\mathscr M$ and $\|T_nf\|_\mathscr M\leq\|f\|_\mathscr M$ by \cite{afp}. Since $T_n$ and $P$ are compact operators, we have that
\[\|C_\psi\|_e\leq\|C_\psi-C_\psi T_n-P(C_\psi-C_\psi T_n)\|.\]
Now choose $f_n\in\mathscr M$ with $\|f_n\|_\mathscr M\leq1$ such that
\[\|C_\psi-C_\psi T_n-P(C_\psi-C_\psi T_n)\|\leq\|(C_\psi-C_\psi T_n)f_n-P(C_\psi-C_\psi T_n)f_n\|_\mathscr M+\varepsilon.\]
By Alaoglou's theorem and the observation preceding the proof of theorem \ref{main}, we may assume that there exists $f\in\mathscr M$ such that $f_n\to f$ uniformly on compact sets. Now
\begin{align}
\label{e3}
\|C_\psi\|_e&\leq\|(C_\psi-C_\psi T_n)f_n-P(C_\psi-C_\psi T_n)f_n\|_\mathscr M+\varepsilon\lesssim\|(C_\psi f_n-C_\psi T_nf_n)''\|_1+\varepsilon\\
&=\iint\limits_{U_R}|(f_n\circ\psi-T_nf_n\circ\psi)''|\, dA+\iint\limits_{L_R}|(f_n\circ\psi-T_nf_n\circ\psi)''|\, dA+\varepsilon\notag\\
&=I_1+I_2+\varepsilon\notag.
\end{align}

To estimate $I_1$, we use the reproducing formula
\[(f_n\circ\psi-T_nf_n\circ\psi)''(z)=-\iint\limits_\mathbbm D\frac1{\overline \alpha}\,(\varphi_\alpha\circ\psi)''(z)(f_n-T_nf_n)''(\alpha)\, dA(\alpha).\]
Then 
\[I_1\leq\left(\sup_{|\alpha|\leq s}\iint\limits_{U_R}|(\varphi_{\alpha}\circ\psi)''|\, dA+\sup_{|\alpha|>s}\|(\varphi_\alpha\circ\psi)''\|_{1}\right)\cdot\iint\limits_\mathbbm D\frac1{|\alpha|}\, |(f_n-T_nf_n)''(\alpha)|\, dA(\alpha).\]
For $|\alpha|\leq s$ we have
\begin{align*}
 \iint\limits_{U_R}|(\varphi_{\alpha}\circ\psi)''|\, dA&=\iint\limits_{U_R}\left|\frac{1-|\alpha|^2}{(1-\overline \alpha\psi)^2}\psi''+2\overline a\frac{1-|\alpha|^2}{(1-\overline \alpha\psi)^3}(\psi')^2\right|\, dA\\
&\lesssim\frac1{(1-s)^3}\iint\limits_{U_R}\left(|\psi''|+|\psi'|^2\right)\, dA.
\end{align*}
Moreover
\begin{align*}
\iint\limits_\mathbbm D\frac1{|\alpha|}\, |(f_n-T_nf_n)''(\alpha)|\, dA(\alpha)\lesssim\|(f_n-T_nf_n)''\|_1\lesssim\|f_n-T_nf_n\|_\mathscr M\leq2.
\end{align*}
Therefore
\begin{equation*}
I_1\lesssim\frac1{(1-s)^3}\iint\limits_{U_R}\left(|\psi''|+|\psi'|^2\right)\, dA+\sup_{|\alpha|>s}\|(\varphi_\alpha\circ\psi)''\|_{1}.
\end{equation*}

To estimate $I_2$, note that
\begin{align*}
I_2&\leq\iint\limits_{L_R}\left(|f_n'\circ\psi-(T_nf_n)'\circ \psi|\cdot|\psi''|+|f_n''\circ\psi-(T_nf_n)''\circ \psi|\cdot|\psi'|^2\right)\, dA\\
&\leq C_{R,R'}\left(\sup_{|z|\leq R'}|f_n(z)-f(z)|+\|T_nf-f\|_\mathscr M\right)\left(\|\psi''\|_{1}+\|\psi'\|^2_{2}\right).
\end{align*} 

So, taking the limits as $n\to\infty$, $R\to1$, $s\to1$, $\varepsilon\to0$ (in the indicated order), (\ref{e3}) gives
\[\|C_\psi\|_e\lesssim \limsup_{|\alpha|\to1}\|(\varphi_\alpha\circ\psi)''\|_{1}.\]
\end{proof}

We remark that the ``standard'' way to prove lower bounds for the essential norm is to choose an appropriate normalized sequence $f_n$ which converges weakly to zero. Then 
\[\|C_\psi\|_e\geq\limsup_n\|C_\psi f_n\|.\]
In our case such a sequence does not exist because the minimal space has the Schur property, that is, weak convergence is equivalent to convergence in the norm. Indeed, $\mathscr M$ is isomorphic to the Bergman space $A^1$ (see, for instance, \cite{zhu}). But $A^1$ is isomorphic to $\ell^1$, \cite{wo}, and the latter space is known to have the Schur property.
It might nevertheless be instructive to give a more direct proof of the fact that $\mathscr M$ and $\ell^1$ are isomorphic. Using standard decomposition arguments one shows that there exist a sequence of points $\alpha_n\in\mathbbm D$, and a sequence of functionals $\Lambda_n\in\mathscr M^*$ such that  the operator  
\[T:\ell^1\to\mathscr M,\quad T\lambda=\sum_n\lambda(n)\,\varphi_{\alpha_n},\] 
is onto, and the operator
\[S:\mathscr M\to\ell^1,\quad Sf=(\Lambda_n(f))_{n=1}^\infty,\]
is an isomorphic embedding. This implies that
\[\ell^1=\ker(T)\oplus S(\mathscr M).\]
Thus, $S(\mathscr M)$ is complemented, therefore, by a theorem of Pelczynski \cite{lt}, it is isomorphic to $\ell^1$. Actually, we don't even need to know that $\mathscr M$ and $\ell^1$ are isomorphic in order to establish the Schur property. Here is a completely independent proof. Suppose toward a contradiction that $f_n$ is a sequence  such that $f_n\overset{w}{\longrightarrow}0$ and $\|f_n''\|_1=1$. Weak convergence implies uniform convergence on compact sets, therefore we can find a subsequence $f_{k_n}$ and an increasing sequence $R_n\in(0,1)$ such that
\[\iint\limits_{|z|<R_{n-1}}|f_{k_n}''(z)|\, dA(z)<\frac1{10},\quad \iint\limits_{R_{n-1}<|z|<R_n}|f_{k_n}''(z)|\, dA(z)>\frac9{10}.\]
Now choose  real functions $\theta_n$ so that
\[|f_{k_n}''(z)|=f_{k_n}''(z)\, e^{i\theta_n(z)},\quad R_{n-1}<|z|<R_n,\]
and put
\[h=\sum_ne^{i\theta_n}\chi_{\{z:R_{n-1}<|z|<R_n\}}.\]
Then $h\in L^\infty(\mathbbm D)$, hence
\[\Lambda_h(f)=\iint\limits_{\mathbbm D}f''h\, dA,\]
is an element of $\mathscr M^*$. 
However
\begin{align*}
\left|\iint\limits_{\mathbbm D}f_{k_n}''h\, dA\, \right|&\geq\left|\,\, \iint\limits_{R_{n-1}<|z|<R_n}f_{k_n}''(z)\, e^{i\theta_n(z)}\, dA(z)\, \right|\\
&-\iint\limits_{|z|<R_{n-1}}|f_{k_n}''(z)|\, dA(z)-\iint\limits_{|z|>R_n}|f_{k_n}''(z)|\, dA(z)\\
&>\frac9{10}-\frac1{10}-\frac1{10},
\end{align*}
contradicting $\Lambda_h(f_{k_n})\to0$.

We further remark that one may use the argument in the proof of the lower bound in theorem \ref{main} to show that if $\psi$ induces a bounded composition operator on the Bergman space $A^1$ then 
\[\|C_\psi\|_e\geq\limsup_{|\alpha|\to1}\|u_\alpha\circ\psi\|_1,\]
where
\[u_\alpha(z)=\frac{1-|\alpha|^2}{(1-\overline\alpha z)^3},\quad \alpha,z\in\mathbbm D.\]
Combining this with the growth estimate
\[|f(z)|\leq\frac{\|f\|_1}{(1-|z|^2)^2}\]
for $f$ in $A^1$, we obtain
\[\|C_\psi\|_e\geq\limsup_{|\alpha|\to1}\left(\frac{1-|\alpha|^2}{1-|\psi(\alpha)|^2}\right)^2.\]
This complements a result due to Vukoti\'c \cite{vukotic}, who proved that in the case of the reflexive Bergman spaces $A^p$, $p>1$, we have
\[\|C_\psi\|_e\geq\limsup_{|\alpha|\to1}\left(\frac{1-|\alpha|^2}{1-|\psi(\alpha)|^2}\right)^{2/p}.\]
Moreover, using the techniques in the proof of the upper bound in  theorem \ref{main} and the decomposition
\[f=\sum_{n=1}^\infty\lambda(n)u_{\alpha_n},\quad\lambda\in\ell^1,\ \alpha_n\in\mathbbm D,\]
for functions in $A^1$, we can prove the corresponding upper bound for composition operators on $A^1$, namely
\[\|C_\psi\|_e\lesssim\limsup_{|\alpha|\to1}\|u_\alpha\circ\psi\|_1.\]
We omit the details.

Returning to the minimal space, we have already mentioned that the  result in \cite{shapiro} implies that if $\|\psi\|_\infty=1$ then $C_\psi$ is not compact. We will now show that something much stronger holds. Namely, the essential norm of $C_\psi$  is at least $1$. We will need the following.
\begin{lemma}
\label{lemma1}
 Let $f\in\mathscr M$ be such that $f(0)=f'(0)=0$ and $f(1)=1$. Then for every $0<r\leq1/2$ we have that
\[\iint\limits_{\Omega_r}|f''|\, dA\gtrsim r,\]
where 
\[\Omega_r=\{z:|z|\leq r\}\cup\bigcup_{|z|\leq r}[z,1].\]
\end{lemma}

\begin{proof}
 An integration by parts yields
\[|f(1-\varepsilon)|-\varepsilon\,|f'(1-\varepsilon)|\leq\int_0^{1-\varepsilon}(1-u)\,|f''(u)|\, du,\]
for small $\varepsilon>0$. Letting $\varepsilon\to0$ and using the fact that $f$ is in the little Bloch space, we get
\begin{equation*}
 \label{e6}
1\leq\int_0^1(1-u)\,|f''(u)|\,du.
\end{equation*}
Therefore
\begin{align*}
 1&\leq\frac 1{\pi r^2}\int_0^1\frac 1{(1-u)}\iint\limits_{|z-u|\leq (1-u)r}|f''(z)|\, dA(z)\, du\\
&\leq\frac 1{\pi r^2}\iint\limits_{\Omega_r}|f''(z)|\int_{|z-1|/(1+r)}^{|z-1|/(1-r)}\frac{dt}t\, dA(z)\lesssim\frac 1r\iint\limits_{\Omega_r}|f''(z)|\, dA(z),
\end{align*}
and we are done.
\end{proof}

\begin{corollary}
\label{cor1}
Let $C_\psi$ be bounded. If $\|\psi\|_\infty=1$ then $\|C_\psi\|_e\geq1$. 
\end{corollary}

\begin{proof}
By automorphism invariance, we may assume that $\psi(0)=0$ and $\psi(1)=1$. Let $\beta_n=1-1/n$ and consider the functions
\[f_n(z)=\frac{g_n((1-\beta_n)z+\beta_n)-(1-\beta_n)g_n'(\beta_n)z}{1-(1-\beta_n)g_n'(\beta_n)},\ |z|<1,\]
where 
\[g_n=\frac{1-\overline{\psi(\beta_n)}}{\psi(\beta_n)-1}\,\varphi_{\psi(\beta_n)}\circ\psi.\]
Then $f_n(0)=f_n'(0)=0$ and $f_n(1)=1$. So, using lemma \ref{lemma1} with, say, $r=1/2$, we get
\begin{align*}
 1&\lesssim\iint\limits_\mathbbm D|f_n''(z)|\, dA(z)=\frac{(1-\beta_n)^2}{|1-(1-\beta_n)g_n'(\beta_n)|}\iint\limits_\mathbbm D|g_n''((1-\beta_n)z+\beta_n)|\, dA(z)\\
&\leq\frac1{|1-(1-\beta_n)g_n'(\beta_n)|}\iint\limits_\mathbbm D|(\varphi_{\psi(\beta_n)}\circ\psi)''(z)|\, dA(z).
\end{align*}
By the Schwarz-Pick lemma
\[|g_n'(\beta_n)|=\frac{|\psi'(\beta_n)|}{1-|\psi(\beta_n)|^2}\leq\frac1{1-\beta_n^2}.\]
Therefore
\[1-\frac{1-\beta_n}{1-\beta_n^2}\lesssim \iint\limits_\mathbbm D|(\varphi_{\psi(\beta_n)}\circ\psi)''(z)|\, dA(z).\]
Letting $n\to+\infty$ we obtain
\[\|C_\psi\|_e\gtrsim1.\]
Note that $\|C_{\psi_n}\|_e\leq\|C_\psi\|_e^n$, where $\psi_n$ is the $n$-fold self-composition $\psi_n=\psi\circ\cdots\circ\psi$. So, we actually have
\[\|C_\psi\|_e\geq1.\]
\end{proof}

This is reminiscent of the situation in $H^\infty$, where the essential norm of a non-compact composition operator is exactly $1$, see \cite{zheng}. In our case, however, $\|C_\psi\|_e$ may take on arbitrarily large values.
\begin{theorem}
\label{blaschke}
 If $B$ is a Blaschke product of degree $n$ then 
$\|C_B\|_e\simeq\|C_B\|\simeq n$.
\end{theorem}

\begin{proof}
 To prove the essential norm estimate, we may assume that $B(0)=0$ because
\[\limsup_{|\alpha|\to1}\|(\varphi_\alpha\circ B)''\|_{1}=\limsup_{|\alpha|\to1}\|(\varphi_\alpha\circ\varphi_\beta\circ B)''\|_{1}\]
for any fixed $\beta\in\mathbbm D$, and moreover $\varphi_\beta\circ B$ is a Blaschke product of degree $n$. Then 
\begin{align*}
 \|(\varphi_\alpha\circ B)''\|_1=\iint\limits_\mathbbm D\left|\,2\overline \alpha\,(1-|\alpha|^2)\, \frac{B'(z)^2}{(1-\overline\alpha B(z))^3}+(1-|\alpha|^2)\, \frac{B''(z)}{(1-\overline\alpha B(z))^2}\right|\, dA(z).
\end{align*}
 Using the change of variable formula for $n$-valent functions and subordination, we estimate
\[(1-|\alpha|^2)\iint\limits_{\mathbbm D}\frac{|B'(z)|^2}{|1-\overline\alpha B(z)|^3}\, dA(z)=n\, (1-|\alpha|^2)\iint\limits_{\mathbbm D}\frac{dA(z)}{|1-\overline\alpha z|^3}\simeq n,\]
\begin{align*}
(1-|\alpha|^2)\iint\limits_{\mathbbm D}\frac{|B''(z)|}{|1-\overline\alpha B(z)|^2}\, dA(z)&\leq(1-|\alpha|^2)\,\|B''\|_\infty\iint\limits_{\mathbbm D}\frac{dA(z)}{|1-\overline\alpha z|^2}\\
&\simeq\|B''\|_\infty(1-|\alpha|^2)\, \log\frac e{1-|\alpha|^2}.
\end{align*}
Therefore
\[\limsup_{|\alpha|\to1}\|(\varphi_\alpha\circ B)''\|_1\simeq n.\]

To prove the norm estimate, we write
\[B=\prod_{j=1}^n\psi_j,\]
where $\psi_j$ is a M\"obius function $\varphi_{\alpha_j}$ times a unimodular constant, and introduce the notation
\[B_j=\prod_{k\ne j}\psi_k.\]
Then 
\[B''=\frac{(B')^2}{B}+\sum_jB_j\,\psi_j''-\sum_jB_j\,\frac{(\psi_j')^2}{\psi_j}.\]
Hence
\[|B''|\leq\frac{|B'|^2}{|B|}+\sum_j|\varphi_{\alpha_j}''|+\sum_j\frac{|\varphi_{\alpha_j}'|^2}{|\varphi_{\alpha_j}|}.\]
Consequently
\[\|B\|_\mathscr M\simeq\iint\limits_\mathbbm D|B''|\, dA+|B(0)|+|B'(0)|\lesssim n.\]
Therefore
\[\sup_B\|C_B\|=\sup_B\sup_{|\alpha|<1}\|\varphi_\alpha\circ B\|_\mathscr M=\sup_B\|B\|_\mathscr M\lesssim n.\]
Since $\|C_B\|_e\leq\|C_B\|$, we see that $\|C_B\|$ is in fact comparable to $n$.
\end{proof}

Note that the argument in the proof of theorem \ref{blaschke} shows that if $\|\psi''\|_\infty<\infty$ then 
\[\|C_\psi\|_e\lesssim\limsup_{|z|\to1}n_\psi(z).\]
 This leads to the conjecture that, at least for ``nice'' symbols, we have 
\[\|C_\psi\|_e\simeq\limsup_{|z|\to1}n_\psi(z).\]
We do not know how to prove (or disprove) this. We will, nevertheless, show that the lower bound in corollary \ref{cor1} can be improved under certain assumptions on the valency of $\psi$ near the boundary. First, we need some preparation. Let $\psi\in\mathscr M$ and \[\Omega(\beta;w;s)=\{z : |z-\beta|\leq s\}\cup\bigcup_{|z-\beta|\leq s}[z,w],\]
 where $\{z : |z-\beta|\leq s\}\subseteq\mathbbm{D}$, $w\in\overline{\mathbbm{D}}$ and $s\leq\frac 12|w-\beta|$. We also let $\alpha=\psi(\beta)$ and we apply lemma \ref{lemma1}, after a trivial normalization, to get
$$\iint\limits_{\Omega(\beta;w;s)}\left|(\phi_{\alpha}\circ\psi)''(z)\right|\,dA(z)\geq c_0s\left|\frac{\alpha-\psi(w)}{(w-\beta)(1-\overline{\alpha}\psi(w))}-\frac{\psi'(\beta)}{|\alpha|^2-1}\right|,$$
Where $c_0$ is an absolute constant. In particular, if $|w-\beta|\geq 1-|\beta|$, then
$$\iint\limits_{\Omega(\beta;w;\frac{1-|\beta|}2)}\left|(\phi_{\alpha}\circ\psi)''(z)\right|\,dA(z)\geq \frac{c_0}2(1-|\beta|)\left|\frac{\alpha-\psi(w)}{(w-\beta)(1-\overline{\alpha}\psi(w))}-\frac{\psi'(\beta)}{|\alpha|^2-1}\right|.$$

\begin{lemma}
\label{lemma2}
Let $\zeta\in\mathbbm{T}$ and suppose that $\psi(\zeta)\in\mathbbm{T}$. Then, for every sequence  $\beta_n$ converging to $\zeta$, if we put $\alpha_n=\psi(\beta_n)$ we have that  either
$$\limsup\frac{|\zeta-\beta_n|}{1-|\beta_n|}\iint\limits_{\Omega(\beta_n;\zeta;\frac{1-|\beta_n|}2)}\big|(\phi_{\alpha_n}\circ\psi)''(z)\big|\,dA(z)\geq\frac{c_0}4$$
or
$$\limsup\frac{|\zeta-\beta_n|}{1-|\beta_n|}\iint\limits_{\Omega(\beta_n;w;\frac{1-|\beta_n|}2)}\big|(\phi_{\alpha_n}\circ\psi)''(z)\big|\,dA(z)\geq\frac{c_0}4$$
for all $w\in\mathbbm{D}$.
\end{lemma}

\begin{proof} To get a contradiction we assume that
$$\limsup\frac{|\zeta-\beta_n|}{1-|\beta_n|}\iint\limits_{\Omega(\beta_n;\zeta;\frac{1-|\beta_n|}2)}\left|(\phi_{\alpha_n}\circ\psi)''(z)\right|\,dA(z)<\frac{c_0}4$$
and
$$\limsup\frac{|\zeta-\beta_n|}{1-|\beta_n|}\iint\limits_{\Omega(\beta_n;w;\frac{1-|\beta_n|}2)}\left|(\phi_{\alpha_n}\circ\psi)''(z)\right|\,dA(z)<\frac{c_0}4$$
for at least one $w\in\mathbbm{D}$.
 These imply
$$\frac 12>\limsup|\zeta-\beta_n|\left|\frac{\alpha_n-\psi(\zeta)}{(\zeta-\beta_n)(1-\overline{\alpha_n}\psi(\zeta))}-\frac{\psi'(\beta_n)}{|\alpha_n|^2-1}\right|$$
and
$$\frac 12>\limsup|\zeta-\beta_n|\left|\frac{\alpha_n-\psi(w)}{(w-\beta_n)(1-\overline{\alpha_n}\psi(w))}-\frac{\psi'(\beta_n)}{|\alpha_n|^2-1}\right|.$$
Adding, we get
$$1>\limsup\left|\frac{\alpha_n-\psi(\zeta)}{1-\overline{\alpha_n}\psi(\zeta)}-\frac{(\zeta-\beta_n)(\alpha_n-\psi(w))}{(w-\beta_n)(1-\overline{\alpha_n}\psi(w))}\right|$$
and, hence,
$$1>1.$$
\end{proof}

Now, let $0<t\leq 1$ and $n\in\mathbbm{N}$. A point $\xi\in\mathbbm{T}$ is called an $(n;t)$ value of $\psi$ if there are $\zeta_1,\ldots,\zeta_n\in\mathbbm{T}$ and a sequence of points $\alpha_m$  converging to $\xi$, such that for each $\alpha_m$  there exist $\beta_{m,1},\ldots,\beta_{m,n}$ so that
$$\psi(\beta_{m,1})=\cdots=\psi(\beta_{m,n})=\alpha_m$$
and 
$$\frac{1-|\beta_{m,j}|}{|\zeta_j-\beta_{m,j}|}\geq t\qquad (1\leq j\leq n).$$
In such a case, it is obvious that $\beta_{m,j}\to\zeta_j\, (1\leq j\leq n)$ and that
$$\psi(\zeta_1)=\cdots=\psi(\zeta_n)=\xi.$$
So, geometrically, $\xi$ can be approximated by points whose preimages lie within a fixed number of Stoltz domains.
\begin{theorem}
If there is at least one $(n;t)$ value of $\psi$, then $\|C_{\psi}\|_e\geq\frac{c_0}4nt$.  
\end{theorem}

\begin{proof} It is obvious by  lemma \ref{lemma2}, that for each $j=1,\ldots,n$ we can choose a fixed $w_j$ so that either $w_j\in\mathbbm{D}$ or $w_j=\zeta_j$ and so that, for some subsequence $\alpha_{m_k}$ (and the corresponding subsequences  $\beta_{m_k,j}\, (1\leq j\leq n)$), all regions $\Omega(\beta_{m_k,j};w_j;\frac{1-|\beta_{m_k,j}|}2)\, (1\leq j\leq n)$ are mutually disjoint and so that  
$$\limsup_{k}\iint\limits_{\Omega(\beta_{m_k,j};w_j;\frac{1-|\beta_{m_k,j}|}2)}\big|(\phi_{\alpha_{m_k}}\circ\psi)''(z)\big|\,dA(z)\geq\frac{c_0}4t.$$
Adding in $j$, we see that
$$\|C_{\psi}\|_e\geq\limsup_{k}\iint\limits_{\mathbbm{D}}\big|(\phi_{\alpha_{m_k}}\circ\psi)''(z)\big|\,dA(z)\geq\frac{c_0}4nt.$$
\end{proof}

\end{document}